\documentclass[11pt]{article}
\usepackage{amssymb}
\title{ On a variant of
 Gagliardo-Nirenberg  inequality deduced from  Hardy\thanks{
{\bf Running head:} Hardy and Gagliardo-Nirenberg inequalities }}
\author{Agnieszka Ka\l amajska
 and Katarzyna Pietruska-Pa\l{}uba\thanks{This work was partially
 done during the second author's stay at the
  Institute of Mathematics of the Polish Academy of Sciences at Warsaw.}}

\date{}

\renewcommand{\it}{\sl}
\renewcommand{\em}{\sl}
plus 6pt minus 12pt

\newcommand{\barint}{
         \rule[.036in]{.12in}{.009in}\kern-.16in
          \displaystyle\int  }
\def\r{{\mathbb{R}}}

\def\rn{{\mathbb{R}^{n}}}

\def\zi{[0,\infty )}

\newtheorem{theo}{\bf Theorem}[section]

\newtheorem{rem}{\bf Remark}[section]

\newtheorem{ex}{\bf Example}[section]

\newcommand{\ds}{\displaystyle}
\newcommand{\ts}{\textstyle}

\newcommand{\vp}{\varphi}

\newcommand{\rp}{\mathbb{R}_+}

\makeatletter \@addtoreset{equation}{section}

\makeatother

\begin{document}
\maketitle \sloppy
\parindent 1em
{\small Institute of Mathematics, University of Warsaw,  ul.\
Banacha 2, 02--097 Warszawa, Poland}

\begin{abstract}\small
We obtain new variants of weighted Gagliardo-Nirenberg
interpolation inequalities in Orlicz spaces, as a consequence of
weighted Hardy-type inequalities. The weights we consider need not
be doubling.
\end{abstract}

\section{Introduction}Gagliardo-Nirenberg interpolation inequalities have already a
long history and several mathematicians investigated their
numerous variants. Their rudiments can be found in the old papers
of Landau (see e.g. \cite{la}), and now they are often identified
with their classical variant
\begin{equation}\label{gn1}
\| \nabla^{(k)} u\|_{L^q(\Omega)}\le C_1\|
u\|_{L^r(\Omega)}^{1-\frac{k}{m}}\|
\nabla^{(m)}u\|^{\frac{k}{m}}_{L^p(\Omega)} + C_2\|
u\|_{L^r(\Omega)},
\end{equation}
where $\Omega\subseteq \rn$ is a domain with sufficiently smooth
boundary, $u:\Omega\rightarrow \r$ belongs to an appropriate
Sobolev space on $\Omega$, $\frac{1}{q}=(1-\frac{k}{m})\frac{1}{r}
+ \frac{k}{m}\frac{1}{p},$ $0<k<m.$ This inequality for
$\Omega=\r$ and $p=q=r=\infty$ was obtained by
 Kolmogorov
\cite{ko} ($C_2=0$ there), whereas  Gagliardo \cite{g} and
Nirenberg \cite{n1} independently proved its extensions to the
form (\ref{gn1}). We refer to the book \cite{hadamard} for an
extensive description of their historical evolution.

Since the Gagliardo-Nirenberg inequalities involve two
differential operators: $\nabla^{(k)}u$ and $\nabla^{(m)} u$, they
are more difficult to analyse than Hardy-type inequalities, which
involve one diferential operator only. This is one of the reasons
why many questions concerning the validity of the
Gagliardo-Nirenberg inequalities remain unsolved so far. For
example,  one asks about Orlicz-space generalizations of
(\ref{gn1}), which can be further extended to Orlicz spaces
$L^M(\mu)$ with a Radon measure $\mu,$ especially when $\mu$ is
nondoubling.

Interest in inequalities in Orlicz-space setting arise from
 linear and nonlinear PDEs, and calculus of variations, which in turn come from mathematical physics.
  See e.g.\  \cite{alber,ball,cia,gm1,gm2,trib2},
where many motivations for investigating  degenerate PDE's in
Orlicz spaces can be found.

The purpose of this paper is to   show that certain variants of
weighted modular Hardy inequalities
\begin{eqnarray}\label{jedd}
\int_\Omega P(|\nabla \varphi| |u|)d\mu \leq K_1 \int_\Omega
P(A|\nabla u|)d\mu + K_2\int_\Omega M(| u|)d\mu,\;\;\;\;u\in
C_0^\infty(\Omega),
\end{eqnarray}
where $\mu(dx)={\rm e}^{-\varphi (x)}{\rm d}x $, $\varphi$  is
locally Lipschitz, imply variants of Gagliardo-Nirenberg
inequalities for modulars:
\begin{eqnarray}
\int_\Omega M(|\nabla u|)d\mu &\leq& L\int_\Omega
P(|\nabla^{(2)}u|)d\mu+ \int_\Omega Q({B}|u|)d\mu,\ \hbox{\rm
where}\ \theta>0,\nonumber \end{eqnarray} and for norms:
\begin{eqnarray} \|\nabla u\|_{L^M(\Omega,\mu)}&\leq& {L}_1
\sqrt{\|\nabla^{(2)}u\|_{L^P(\Omega,\mu)} \|u\|_{L^Q(\Omega,\mu)}}
+ {L}_2\|u\|_{L^Q(\Omega,\mu)},\label{obiad1}
\end{eqnarray}
valid with general ($u$-independent) constants
$A,K_1,K_2,L,L_1,L_2,B$. Our approach requires $M$ to be an
$N$-function satisfying the $\Delta_2$-condition. Functions $P$
and  $Q$ are  tied with $M$ by  Young-type inequality:
$$\frac{M(u)}{u^2}\,vw\leq M(u)+P(v)+Q(w).$$ For details, see
Theorems \ref{theoA}, \ref{theoAB}, and \ref{theoB}.

As opposed to previous woks of  Gutierrez and Wheeden \cite{gw},
and also Chua  \cite{skch,skch94}, Bang and coauthors
(\cite{bang}-\cite{hama}), the measure $\mu $ considered here
needs not to be doubling. This allows for obtaining inequalities
e.g. for measures with finite mass on unbounded domains, which
were formerly excluded from investigation. In papers
\cite{gw,skch,skch94}, Gagliardo-Nirenberg inequalities were
deduced from {\em local} Poincar\'e inequalities. In present work,
we work with global inequalities only, and show how {\em global}
Hardy-type inequalities result in Gagliardo-Nireberg inequalities.

Inequality (\ref{obiad1}), obtained here as a consequence of Hardy
inequality, extends our former results from \cite{BBMS}. In that
paper, nondoubling measures were considered as well, but
 the conditions for $M,P,Q$ were different.  In particular, the case
 $M=P=Q$ was not permitted  in \cite{BBMS},  see Remark 4.4 in \cite{BBMS}.
  This situation is rectified  by present approach.

 It is our intention to focus on Hardy-type inequalities
 (\ref{jedd}). They imply a big range of Gagliardo-Nirenberg
 inequalities. Since (\ref{jedd}) is  valid for a vast class of
 admissible measures, possibly nondoubling, our approach yields
 Gagliardo-Nirenberg inequalities which are often new also in the
 $L^p-$setting.

For other results concerning the Gagliardo-Nirenberg inequalities
in Orlicz spaces, we refer to \cite{bang}-\cite{hama},
\cite{AKKPPlms}-\cite{koche}.

\section{Preliminaries}

\noindent{\bf Notation}
\newline
In the sequel we assume that $\Omega\subset\rn$ is an open domain.
By $C_0^\infty (\Omega)$ we denote smooth functions compactly
supported in  $\Omega$, we use the standard notation
$W^{m,p}(\Omega)$ and $W^{m,p}_{loc}(\Omega)$ for global and local
variants of  Sobolev spaces. Lower-case symbol $c$ denotes a
universal constant whose value is irrelevant. For important
constants we use the upper-case letters.

\bigskip

\noindent{\bf Orlicz spaces}
\newline
Let us report some basic information about Orlicz spaces,
referring e.g. \cite{kr,rao-ren} for details.

 Suppose that $\mu$ is a positive Radon
measure on $\rp$ and let $M:[0,\infty)\to[0,\infty)$ be an
$N-$function, i.e. a continuous convex function satisfying
$\lim_{\lambda\to 0}\frac{M(\lambda)}{\lambda}=
 \lim_{\lambda\to\infty}\frac{\lambda}{M(\lambda)}=0.$
 The symbol $M^*$ denotes the  function complementary to the given
$N-$function  $M$, i.e. the function $M^*(y)=\sup_{x>0}[xy-M(x)]$,
where $y\geq 0$.
  It is again an $N-$function and
it satisfies the Young inequality:
\[
xy\leq M(x)+M^*(y), \;\;\mbox{ for } \;\; x,y\geq 0.
\]
Given two functions $M_1$ and $M_2$, we write $M_1\asymp M_2$ if
there exist two constants $c_1,c_2$ such that $c_1M_2(\lambda)\le
M_1(\lambda)\le c_2M_2(\lambda)$, for every $\lambda>0$ (or for
every $\lambda$ from indicated range).

Let $\mu$ be a nonnegative Borel measure on $\Omega$.
 The weighted Orlicz space $L^{M}(\Omega,\mu)$ we deal
with is by definition the space
\[L^{M}(\Omega,\mu)\stackrel{def}{=}\{f:\Omega\rightarrow\r\mbox{ measurable} :   \int_{\Omega}
 M(\frac{|f(x)|}{K})d\mu(x)\leq 1\ \mbox{ for some }\ K >0\} ,  \]
 equipped with the Luxemburg norm
\[\|f\|_{L^{M}(\Omega,\mu)}=\inf\{ K>0 : \int_{\Omega}M(\frac{|f(x)|}{K})d\mu(x)\leq
1\} .\] It is a Banach space. When $M(\lambda)=\lambda^p,$ then we
have $L^{M}(\Omega,\mu)= L^p_\mu (\Omega)$ -- the classical $L^p$
space.

The function  $M$ is said to fulfill the $\Delta_2-$condition if
and only if for some constant $c>0$ and every $\lambda >0,$ we
have
\begin{equation}\label{delta2}
M(2\lambda)\leq cM(\lambda).
\end{equation}
In the class of differentiable convex functions the
$\Delta_2-$condition is equivalent to:
\begin{equation}\label{pierwsze}
\lambda M'(\lambda)\leq D M(\lambda),
\end{equation}
satisfied for every $\lambda >0$, with the constant $D$ being
independent of $\lambda$ (see e.g. \cite{kr}, Theorem 4.1). The
smallest possible  constant in (\ref{pierwsze}) will be denoted by
$D_M.$ One considers also the condition
\begin{equation}\label{drugie}
{d} M(\lambda)\leq \lambda M'(\lambda).
\end{equation}
It holds true with some $d>1$  when the dual function, $M^*,$
satisfies the $\Delta_2-$condition  (see e.g. \cite{kr}, Theorem
4.3 or \cite{akikppstud}, Proposition 4.1). The biggest  possible
constant in (\ref{drugie}) will be denoted by $d_M.$ Constants
$d_M$ and $D_M$ are called Simonenko lower and upper indices of
$M$ (see e.g. \cite{boy}, \cite{simon}).

When $M$ is such an $N-$function that both  (\ref{pierwsze})
 and (\ref{drugie}) are satisfied, then
\begin{eqnarray}
 M(a\lambda)&\leq& {\rm max}(a^{d_M}, a^{D_M})M(\lambda)=: \bar c(a)M(\lambda),
\label{deltadwa}
\end{eqnarray}
 for every  $\lambda>0, a>0$ (see e.g. \cite{akikppstudia}, Lemma
4.1, part iii)).

 We will need the following property of modulars:
 \begin{equation}\label{norm1}
 \int_{\Omega} M(\ts\frac{f(x)}{\|f\|_{L^M(\Omega,\mu)}})\,d\mu(x)\leq
 1.
 \end{equation}
 When $M$ satisfies the $\Delta_2-$condition, then
 (\ref{norm1}) becomes an equality.

\noindent {\bf The  assumptions}

The assumptions considered in the sequel are as follows.

\begin{description}
\item[(M)]
$M:\zi\to\zi$ is an $N-$function of class $C^1((0,\infty))$,
satisfying the $\Delta_2-$condition and such that $\frac{M'(
\lambda)}{{\lambda}}$ is bounded next to zero;
\item[($\mu$)]
$\mu(dx)={\rm exp}(-\varphi(x))dx$ is a Radon measure on $\Omega,$
where $\varphi:\Omega\to\mathbb{R}$ belongs to
$W^{1,\infty}_{loc}(\Omega)$;
\item[(Y)] $P$ and $Q$ are two real nonnegative and nondecreasing
measurable functions on $\zi$ with $P(0)=Q(0)=0,$  such that for
any $u,v,w>0$ the following Young-type inequality holds
\begin{equation}\label{assumA1}
\frac{M(u)}{u^2}\,vw\leq M(u)+P(v)+Q(w).
\end{equation}
\end{description}

The inequality  {\bf (Y)} is fulfilled for example when the
following condition holds (see  \cite{AKKPPIndiana},
 Cor. 4.1):
\begin{description}
\item [(MF)] $M$ is an $N-$function satisfying the $\Delta_2-$condition and such that
${M(\lambda)}/{\lambda^2}$ is nondecreasing, $ P(\lambda)=
CM({F(\sqrt\lambda)}),$ $Q(\lambda)= CM({F^*(\sqrt\lambda)}),$
where $F$ is another $N-$function.
 \end{description}

\begin{rem}\label{konsul}\rm
\begin{enumerate}
\item The choice of $F(\lambda)=\frac{\lambda^2}{2}$  in {\bf
(MF)} gives that {\bf (Y)} holds   with $P=Q=M$  (with $C=1$).
\item
 Suppose that $M(\lambda)=\lambda^p,$ $p\geq 2.$ Choose $F(\lambda)=
\frac{\lambda^s}{s}$ with $s=\frac{2q}{p}$
 to obtain  $P(\lambda)=C\lambda^q$,
$Q(\lambda)=C\lambda^r$, where
$\frac{2}{p}=\frac{1}{q}+\frac{1}{r}$
 (the classical Gagliardo-Nirenberg triple),
\item When   $M(\lambda)= \lambda^p\ln(2+\lambda)^\alpha, p\geq 2,
\alpha>0,$ then the choice of
$F(\lambda)=\lambda^s\ln(2+\lambda)^\mu,$ with $s=\frac{2q}{p},$
$\mu=\frac{\beta-\alpha}{p}$ results in $P(\lambda)\asymp
\lambda^q\ln(2+\lambda)^\beta, $
 $Q(\lambda)\asymp \lambda^r\ln(2+\lambda)^\gamma,$
where the parameters are related through
$\frac{2}{p}=\frac{1}{q}+\frac{1}{r},$
$\frac{2\alpha}{p}=\frac{\beta}{q}+ \frac{\gamma}{r}$ (the
logarithmic
  Gagliardo-Nirenberg triple considered in \cite{AKKPPlms} and \cite{AKKPPcollmath}).
\item Inequality (\ref{assumA1}) can be obtained from
the multiversion of the Young inequality due to Cooper:
\[
\Pi_\nu a_\nu \le\sum_{\nu} \int_0^{a_\nu}\frac{g_\nu (x)}{x}dx,
\]
where $\Pi_{\nu}g_{\nu}^{-1}(x)=x$, $a_\nu\ge 0$ and $g_\nu$ are
continuous strictly increasing functions such that $g_i(\zi)=\zi$
(see Theorem 159 in \cite{hlp} and Theorem 4.3 in
\cite{AKKPPIndiana}).
\end{enumerate}
\end{rem}

\section{Main results}

Our goal now is to show that certain  Hardy-type inequalities
imply the Gagliardo-Nirenberg ones. Let us start with the
following result.

\begin{theo}\label{theoA} Let $M,P,Q$ be three functions on $\zi$
satisfying {\bf (M), (Y)}, and let $\mu$ be a Radon measure on
$\zi$ satisfying {\bf ($\mu$)}. Suppose that the the following
Hardy-type inequality holds true:
\begin{description}
\item[(H)]
for any $u\in C_0^\infty(\Omega)$
\begin{equation}\label{assumA2}
\int_\Omega P(|\nabla \varphi| |u|)d\mu \leq K \int_\Omega
P(A|\nabla u|)d\mu
\end{equation}
with  positive constants $K,A$ not depending on $u.$
\end{description}
 Then we have:
\begin{description}
\item[1)] there exist constants $L,B>0$ such that for any  $\theta >0$ and  any $u\in
C_0^\infty(\Omega)$
\begin{equation}\label{statA1}
\int_\Omega M(|\nabla u|)d\mu \leq L\int_\Omega
P(\theta|\nabla^{(2)}u|)d\mu+ \int_\Omega
Q(\frac{B}{\theta}|u|)d\mu,
\end{equation}
\item[2)] if additionally  $P$ and $Q$ are $N-$functions, then for any $u\in
C^\infty_0(\Omega)$ we have
\begin{equation}\label{statA3}
  \|\nabla u\|_{L^M(\Omega,\mu)}\leq \tilde{L}
\sqrt{\|\nabla^{(2)}u\|_{L^P(\Omega,\mu)}
\|u\|_{L^Q(\Omega,\mu)}},
\end{equation}
where $\tilde{L}=2(L+2)\sqrt{B}$, $L$ and $B$ are the same as in
(\ref{statA1}).
\end{description}
\end{theo}

\begin{rem}\rm Under the assumptions of Theorem \ref{theoA}, when either $P$ or $Q$  satisfies the
$\Delta_2-$condition, then we have
\begin{equation}\label{statA4}
\int_\Omega M(|\nabla u|)d\mu \leq  L_1\int_\Omega
P(|\nabla^{(2)}u|)d\mu+ L_2 \int_\Omega
 Q(|u|)d\mu,
\end{equation}
with $L_1,L_2$ independent of $u.$
\end{rem}

\noindent{\bf Proof.} We start with the following inequality
(Lemma 3.1 of \cite{BBMS}), valid for any $u\in
C_0^\infty(\Omega)$:
 \begin{eqnarray}
  {I}&=&\int_{\Omega} M(|\nabla u|)d\mu\label{gauss1}\\
& \leq & \alpha_n\int_{\Omega\cap \{ \nabla u\neq 0 \}
}\frac{M(|\nabla u|)}{|\nabla u|^2}|\nabla^{(2)}u|\,|u| \,d\mu  +
\int_{\Omega\cap \{ \nabla u\neq 0 \} }\frac{M(|\nabla
u|)}{|\nabla u|}|\nabla\vp|\,|u|d\mu \nonumber\\[2mm]
{} & =: &\alpha_n{ I}_1+{ I}_2,\nonumber
\end{eqnarray}
where $\alpha_n$ depends on $n$ only. In \cite{BBMS}, the proof is
given for $\Omega=\mathbb{R}^n,$ and $\varphi\in C^1(\rn).$ It
requires only minor alterations to cover the present case.

To estimate $I_1$ and $I_2,$ we use the assumption (\ref{assumA1})
twice. One has, for any given $\epsilon,\theta>0:$

\begin{eqnarray}
I_1&=& \epsilon \int_{\Omega\cap \{ \nabla u\neq 0 \}
}\left(\frac{M(|\nabla u|)}{|\nabla u|^2}\right)\left( \theta
|\nabla^{(2)}u|\right)\,\left(\frac{|u|}{\theta\epsilon}\right)
\,d\mu \noindent\nonumber\\
 &\leq& \epsilon  \int_{\Omega} M (|\nabla u|)d\mu
+ \epsilon\int_{\Omega} P (\theta| \nabla^{(2)} u|)d\mu  +
\epsilon\int_{\Omega} Q (\frac{| u|}{\theta\epsilon})\,d\mu,
\label{pott}
\end{eqnarray}
and
\begin{eqnarray}
{I}_2&=& A\epsilon \int_{\Omega\cap \{ \nabla u\neq 0 \} }\left(
\frac{M(|\nabla u|)}{|\nabla u|^2}\right)\left( |
\nabla\vp|\frac{\theta}{A}| \nabla
u|\right)\left(\frac{|u|}{\theta\epsilon}\right)d\mu
\le\nonumber\\
&\le& A\epsilon  \int_{\Omega} M (|\nabla u|)d\mu + A\epsilon
\int_{\Omega} P ( |\nabla\vp|\frac{\theta}{A}|\nabla  u|)d\mu
+A\epsilon \int_{\Omega} Q (\frac{|
u|}{\theta\epsilon})\,d\mu.\label{potrz}
\end{eqnarray}
To estimate the central term in (\ref{potrz}), we apply inequality
(\ref{assumA2}) to the function $\frac{\theta}{A} f(x)$ where
$f(x)=|\nabla u(x)|.$ Whenever $\nabla u\neq 0,$ one has:
\[\frac{\partial}{\partial x_i}f(x)= \langle \frac{\nabla u(x)}{|\nabla u(x)|},
\nabla \frac{\partial}{\partial x_i}u(x)\rangle,\] and so
\[|\nabla f(x)|^2=\sum_{i=1}^n |\langle \frac{\nabla u(x)}{|\nabla
u(x)|}, \frac{\partial}{\partial x_i}(\nabla u(x))\rangle|^2\leq
\sum_{i=1}^n  \|\frac{\partial}{\partial x_i}(\nabla u(x))\|^2
=\|\nabla^{(2)}u(x)\|^2.
\]
 Since $P$ is nondecreasing, we get from {\bf (H)}:
\[ \int_\Omega
P(|\nabla\varphi|\, \frac{\theta}{A} |\nabla  u|)d\mu\leq
K\int_\Omega P(\theta |\nabla^{(2)}u|)d\mu.\] Using this fact,
summing up estimates (\ref{pott}) and (\ref{potrz}) we obtain:
\[
I\leq \epsilon(\alpha_n+A)I + \epsilon(\alpha_n+KA)\int_\Omega
P(\theta|\nabla^{(2)}u|)\,d\mu+ \epsilon(\alpha_n+A)\int_\Omega
Q(\frac{|u|}{\theta \epsilon})\,d\mu.
\]
Choose $\epsilon=\frac{1}{2(\alpha_n+A)},$ which after rearranging
gives
\[I\leq (K+1) \int_\Omega P(\theta |\nabla^{(2)}u|)\,d\mu + \int_\Omega Q\left(\frac{2(\alpha_n+A)}{\theta}|u|\right)d\mu.\]
This proves statement {\bf 1)}.

\smallskip
\noindent To prove statement {\bf 2)} we take
 an arbitrary $u\in
C_0^\infty(\Omega)$ and apply (\ref{statA1}) to \[ \widetilde
u:=\frac{u}{a+b},\ \ {\rm where} \ \
a:={\|\theta\nabla^{(2)}u\|_{L^P(\Omega,\mu)},\ \
b:=\|\frac{B}{\theta}u(x)\|_{L^Q(\Omega,\mu)}}.\] Without loss of
generality we can assume that both numbers $a$ and $b$ are not
zero, because when either $a$ or $b$ equals to zero, then we must
have $\nabla u(x)=0$ a.e. and (\ref{statA3}) follows trivially.
Inequality (\ref{statA1}) for $\widetilde u$ reads
\begin{eqnarray*}
\int_\Omega M(|\nabla \widetilde u|)d\mu & \leq &  L\int_\Omega
P(\ts\frac{|\theta\nabla^{(2)}u|}{a+b})d\mu + \ds\int_\Omega
Q(\ts\frac{|\frac{B}{\theta}u|}{a+b})d\mu\\ & \le & L\int_\Omega
P(\ts\frac{|\theta\nabla^{(2)}u|}{a})d\mu + \ds\int_\Omega
Q(\ts\frac{|\frac{B}{\theta}u|}{b})d\mu
\\ &=& L\int_\Omega
P(\ts\frac{|\theta\nabla^{(2)}u|}{\|\theta\nabla^{(2)}u\|_{L^P(\Omega,\mu)}})d\mu
+ \ds\int Q(\ts\frac{|\frac{B}{\theta}u|}{\| (B/\theta)
u\|_{L^Q(\Omega,\mu)}})d\mu
 \leq L+1.
\end{eqnarray*}
In the last inequality  we have used the property \( \int_\Omega
R(\ts\frac{w}{\|w\|_{L^R(\Omega,\mu)}})d\mu\le 1 \) of modular
functionals. Since for any $f\in L^M(\Omega,\mu)$ one has
$\|f\|_{L^M(\Omega,\mu)}\leq \int_\Omega M(|f|)d\mu +1$ (see (9.4)
and (9.20) of \cite{kr}), this gives
\[\|\nabla \tilde u\|_{L^M(\Omega,\mu)}\leq L+2.\]  Consequently,
\[
\|\nabla  u\|_{L^M(\Omega,\mu)}\leq (L+2)(a+b)= (L+2)\left( \theta
\|\nabla^{(2)}u\|_{L^P(\Omega,\mu)} +
\frac{B}{\theta}\|u(x)\|_{L^Q(\Omega,\mu)}\right).
\]
Minimizing the right hand side with respect to $\theta$ gives the
result. The proof is complete.\hfill $\Box$

Our next theorem covers the case when the Hardy inequality {\bf
(H)} does not hold, but it holds when an extra term, depending on
$u$, is added to the right-hand side.

\begin{theo}\label{theoAB}
Suppose that   the assumptions of Theorem \ref{theoA} are
satisfied, and the following Hardy-type inequality holds true:
\begin{description}
\item[(H1)] for
any $u\in C_0^\infty(\Omega)$,
\begin{equation}\label{assumA23}
\int_\Omega P(|\nabla \varphi| |u|)d\mu \leq K_1 \int_\Omega
P(A|\nabla u|)d\mu + K_2\int_\Omega M(| u|)d\mu
\end{equation}
\end{description}
with  positive constants $K_1,K_2,A$ not depending on $u.$
\newline
 Then we have:
\begin{description}
\item[1)] there exist constants $L,B>0$ such that for any  $\theta \in (0,1]$ and any $u\in
C_0^\infty(\Omega)$
\begin{equation}\label{statA13}
\int_\Omega M(|\nabla u|)d\mu \leq L\int_\Omega
P(\theta|\nabla^{(2)}u|)d\mu+ \int_\Omega
Q(\frac{B}{\theta}|u|)d\mu,
\end{equation}
\item[2)] if $P$ and $Q$ are $N-$functions, then for any $u\in
C^\infty_0(\Omega)$ we have
\begin{equation}\label{statA33}
\|\nabla u\|_{L^M(\Omega,\mu)}\leq {L}_1
\sqrt{\|\nabla^{(2)}u\|_{L^P(\Omega,\mu)} \|u\|_{L^Q(\Omega,\mu)}}
+ {L}_2\|u\|_{L^Q(\Omega,\mu)},
\end{equation}
where ${L}_1=2(L+2)\sqrt{B}$,  ${L}_2=2(L+2){B}$, $L$ and $B$ are
the same as in (\ref{statA13}).
\end{description}
\end{theo}

\noindent{\bf Proof.}\\
\phantom{}\  {\bf 1)} We start with inequality (\ref{gauss1}) and
repeat the arguments in the proof of Theorem \ref{theoA} up to the
formula (\ref{potrz}). Now, instead of (\ref{assumA2}), we apply
(\ref{assumA23}) to the function $f(x)=\frac{\theta}{A}|\nabla
u(x)|.$ Since we have assumed $\theta\leq 1,$ we  get:
\[
\int_\Omega P(|\nabla\varphi|\, \frac{\theta}{A} |\nabla
u|)d\mu\leq K_1\int_\Omega P(\theta |\nabla^{(2)}u|)d\mu +
\bar{K}_2\int_{\Omega} M(|\nabla u|)d\mu,
\]
where $\bar{K}_2 =\bar{c}(\frac{1}{A})K_2$, $\bar{c}(\cdot )$
comes from (\ref{deltadwa}). This, (\ref{pott}) and (\ref{potrz})
 lead to the inequality
\begin{eqnarray*}
I\leq
 \epsilon(\alpha_n+A +\bar{K}_2A)I
+ \epsilon(\alpha_n+K_1A)\int_\Omega
P(\theta|\nabla^{(2)}u|)\,d\mu+ \epsilon(\alpha_n+A)\int_\Omega
Q(\frac{|u|}{\theta \epsilon})\,d\mu.
\end{eqnarray*}
The choice of $\epsilon = \left( {2(\alpha_n+A +
A\bar{K}_2)}\right)^{-1}$ implies
\begin{equation}\label{trzecia}
I\le L\int_\Omega P(\theta |\nabla^{(2)}u|)d\mu + \int_\Omega
Q\left(\frac{B}{\theta}|u|\right)d\mu ,
\end{equation}
where $L=K_1+1$, $B=2(\alpha_n +A +A\bar{c}(\frac{1}{A}){K}_2)$.
This completes the proof of part {\bf 1)}.

{\bf 2)} To prove the second part we observe that arguments
similar to those
  in the proof of second part in Theorem \ref{theoA}
lead to the inequality
\[
\|\nabla u\|_{L^M(\Omega,\mu)} \le\tilde{L}\left( \theta
\|\nabla^{(2)} u\|_{L^P(\Omega,\mu)} +\frac{B}{\theta}\|
u\|_{L^Q(\Omega,\mu)}  \right),
\]
where  $\tilde{L} = L+2$, $L=K_1+1$ is the constant from
(\ref{trzecia}),  $\theta\in (0,1]$ is arbitrary. Minimization of
the inequality $a\le \tilde{L}(\theta b + \frac{1}{\theta}c)$ with
respect to $\theta\in (0,1]$ gives the desired result. Indeed,
when $\bar{\theta}:=\sqrt{\frac{c}{b}}\in (0,1)$, then we get
$a\le 2\tilde{L}\sqrt{bc}$, while in the remaining case $c\ge b$
 we have $a\le 2\tilde{L}c$ (choose ${\theta}=1$). In either case
inequality $a\le 2\tilde{L}(\sqrt{bc}+c)$ holds. This completes
the proof of the theorem. \hfill$\Box$

Since the condition {\bf (Y)} is fulfilled for $P=Q=M$ (see Remark
\ref{konsul}), as an immediate consequence we obtain the
following:

\begin{theo}\label{theoB}
Suppose that $M$ is an $N-$function satisfying condition {\bf
(M)}, and let $\mu$ be a Radon measure on $\zi $ satisfying {\bf
($\mu$)}. Moreover, assume that ${M(\lambda)}/{\lambda^2}$ is
nondecreasing. If for every $u\in C^\infty_0(\Omega)$ the
following
 Hardy-type
inequality holds true:
\begin{equation}\label{assumB1}
\int_\Omega M(|\nabla\varphi|\,|u|)\,d\mu\leq K_1\int_\Omega
M(|\nabla u|)\,d\mu +K_2\int_\Omega M(|u|)\,d\mu,
\end{equation}
then we have:
\begin{description}
\item[1)] there exist positive constants $C_1,C_2$ such that
for  any $\theta\in (0,1]$ and  any $u\in C_0^\infty(\Omega)$
 \begin{equation}\label{statB1} \int_\Omega M(|\nabla u|)\,d\mu\leq
C_1\int_\Omega M(\theta |\nabla^{(2)}u|)\,d\mu +C_2 \int_\Omega
M(\frac{|u|}{\theta})\,d\mu;
\end{equation}
\item[2)] there exist positive constants $\tilde{C}_1,\tilde{C}_2$ such that
for any $u\in C_0^\infty(\Omega)$
 \begin{equation}\label{statB2} \|\nabla
u\|_{L^M(\Omega,\mu)}\leq\tilde C_1
\sqrt{\|\nabla^{(2)}u\|_{L^M(\Omega,\mu)} \|u\|_{L^M(\Omega,\mu)}}
+\tilde C_2 \|u\|_{L^M(\Omega,\mu)}.
\end{equation}
\end{description}
\end{theo}

\begin{rem}\rm
Constants $L,B$ in the inequality (\ref{statA1}) are of the form:
$L=K+1$, $B=2(\alpha_n+A)$, where $K,A$ are the same as in
(\ref{assumA2}), $\alpha_n$ depends on the dimension (it appears
in (\ref{gauss1}) the proof of Theorem \ref{theoA}). It is known
that $\alpha_n\le c\sqrt{n}$, with $c>0$ independent of $n$ (see
Lemma 3.1 in \cite{BBMS}). This is readily seen from the proof of
Theorem \ref{theoA}.
\end{rem}

\begin{rem}\rm
It can be deduced from the proof of Theorem \ref{theoAB} that
constants $L,B$ in the inequality (\ref{statA13}) are of the form:
$L=K_1+1$, $B=2(\alpha_n+A + A\bar{c}(\frac{1}{A})K_2)$, where $A,
K_1, K_2$ are the same as in (\ref{assumA23}), $\bar{c}$ is the
same as in (\ref{deltadwa}),  $\alpha_n$ appears in (\ref{gauss1})
is such that $\alpha_n\le c\sqrt{n}$, with $c>0$ independent of
$n$.
\end{rem}

\begin{rem}[open question]\rm
Condition $M\in\Delta_2$ in the assumptions of Theorem \ref{theoA}
is only needed  to derive the formula (\ref{gauss1}). We do not
know whether one can extend
 (\ref{statA1}), (\ref{statA3}) to functions $M$ for which the
 $\Delta_2-$conditions is not satisfied.
  Some
results concerning the Gagliardo-Nirenberg inequalities
(\ref{statA1}), (\ref{statA3}) hold true without the
$\Delta_2-$condition imposed on $M$, but  for a restricted family
of measures, see e.g. \cite{bang}-\cite{hama},\cite{akikppstud}.
\end{rem}

\bigskip

\section{Discussion and examples}\label{beczka}
Three theorems  from the previous section reduce the question
about the validity of  the Gagliardo-Nirenberg inequality  for
given $N-$functions and measures to a question about the validity
of Hardy-type inequalities. We will discuss it now.

\subsection{The scope of Theorem \ref{theoA}}

\subsubsection{The case $\Omega=\r_+$, condition {\bf (H)}}

A necessary and sufficient  condition for Radon  measures $\mu,
\nu$ to obey the inequality
\begin{equation}\label{inek1}
\int_{\r_+}\left|\int_0^x f(t)dt\right|^p d\nu (x) \le
C\int_{\r_+} |f(x)|^p d\mu (x),
\end{equation}
was given by Muckenhoupt (see \cite{muc} or \cite{ma}, Section
1.3, Theorem 1):
\begin{equation}\label{ineq2}
\sup_{r>0} \left(\nu (r,\infty)\right)^{\frac{1}{p}}\left(\int_0^r
\left(\frac{d\mu^*}{dx}\right)^{-\frac{1}{p-1}}dx\right)^{(p-1)/p}<\infty,
\end{equation}
where $\mu^*$ is the absolutely continuous part of $\mu.$ Since
for $u\in C_0^\infty (0,\infty)$ one has $u(x)=\int_0^x u'(t)dt,$
it follows that whenever $\nu,\mu$ obey (\ref{ineq2}), then the
inequality
\begin{equation}\label{n1}
\int_{\r_+} |u(x)|^pd\nu (x)\le C\int_{\r_+} |u'(x)|^pd\mu (x)
\end{equation}
holds for $u\in C_0^\infty (\r_+)$. Observe that  in the
particular case of $d\nu (x)=|\varphi'(x)|^p {\rm exp}(-\varphi
(x))dx$, $d\mu (x) = {\rm exp}(-\varphi (x))dx,$ inequality
(\ref{n1}) is nothing but our condition {\bf (H)} for
$P(\lambda)=\lambda^p.$  In this case,  condition (\ref{ineq2})
reads:
\begin{equation}\label{ineq3}
\sup_{r>0} \left( \int_r^\infty |\varphi'(x)|^p {\rm exp}(-\varphi
(x))dx \right)^{\frac{1}{p}}\left(  \int_0^r {\rm exp}(-\varphi
(x))^{-\frac{1}{p-1}}dx \right)^{ 1-\frac{1}{p}} <\infty,
\end{equation}
and so when (\ref{ineq3}) holds true, then {\bf(H)} is true for
 $P(\lambda )=\lambda^p$, $\Omega
=\r_+$. Therefore we obtain:
\begin{theo}\label{Mukh}
Let $p>1$ be given, and let $\mu(dx)= e^{-\vp(x)}dx$ be a Radon
measure on $\zi$ satisfying {\bf ($\mu$)}. Suppose that
(\ref{ineq3}) holds true. Next, let $M$ be an $N-$function
satisfying condition {\bf (M)}, and let $Q$ be another
$N-$function, such that {\bf (Y)} is satisfied for $M, $
$P=P(\lambda)=\lambda^p,$ and $Q.$ Then for any $u\in
C_0^\infty(\rp)$ one has
\begin{equation}\label{iii}
\int_{\rp} M(|u'|)d\mu(x)\leq K_1 \int_{\rp} |u''|^p d\mu
+K_2\int_{\rp} Q(|u|)d\mu,
\end{equation}
and also
\begin{equation}\label{iiin}
\|u'\|^2_{L^M(\rp,\mu)}\leq
K\|u''\|_{L^p(\rp,\mu)}\|u\|_{L^Q(\rp,\mu)},
\end{equation}
where the constants $K_1,K_2,K$ do not depend on $u.$
\end{theo}

We illustrate this case with two examples.

\begin{ex}\rm{\bf[classical Hardy inequality]}\\[1mm]
 Consider the classical Hardy
inequality, which deals with power weights \cite{har20},
\cite{hlp}:
\begin{equation}\label{clasha}
\int_{\r_+} |u(t)|^pt^{\alpha -p}dt \le C\int_{\r_+}
|u'(t)|^pt^{\alpha}dt,
\end{equation}
where $C=\left( \frac{p}{|\alpha -p +1| }\right)^p$, $\alpha\neq
p-1$. In this case we have $\mu(dt)= {\rm exp}(-\varphi(t))dt$,
where $\varphi(t)=-\alpha\ln t$. In particular
$\varphi'(t)=-\frac{\alpha}{t}$, and (\ref{clasha}) reads
\[
\int_{\r_+}( |\varphi'(t)|\,|u(t)|)^pt^{\alpha}dt \le
\tilde{C}\int_{\r_+} |u'(t)|^pt^{\alpha}dt,
\]
where $\tilde{C}=\left( \frac{p\alpha}{|\alpha -p +1| }\right)^p$.
\end{ex}

\begin{ex}\rm{\bf [Hardy inequality and power-exponential weights]}\\[1mm]
We now consider  measures on $(0,\infty)$ with
power-exponential-type densities:
\begin{equation}\label{expy}
\mu(dx)= x^\alpha e^{-x^\beta} \,dx={\rm exp}(-\varphi (x))dx,
\;\;\;\;\; \alpha\geq 0,\beta>0, \;\;\;\;\; \varphi (x)=-\ln
x^{\alpha} +x^\beta.
\end{equation}
This class of measures contains in particular the exponential
distribution ($\alpha=0,$ $\beta=1$) and the Gaussian distribution
($\alpha=0,$ $\beta=2$).

As $|\varphi'(x)|\asymp \frac{1}{x}$ for $x$ being small and
$|\varphi'(x)|\asymp {x}^{\beta -1}$ for $x$ close to $\infty$,
the Muckenhoupt condition (\ref{ineq3}) for the measure
(\ref{expy}) is equivalent to:
\begin{eqnarray}
&&\sup_{r>0}
\left(A(r)\right)^\frac{1}{p}\left(B(r)\right)^{1-\frac{1}{p}}<\infty ,\ {\rm where}\ \nonumber\\
&A(r):=& \left(\int_r^1
x^{-p+\alpha}e^{-x^\beta}dx\right)\chi_{r<1} + \left(\int_r^\infty
x^{(\beta-1)p+\alpha}{\rm
e}^{-{x^\beta}}\,dx\right)\chi_{r\ge 1},\nonumber\\
&B(r):=& \int_0^r \frac{ {\rm e}^{\frac{x^\beta}{p-1}}
}{x^{\frac{\alpha}{p-1}}}\,dx.\label{mucexp}
\end{eqnarray}
We observe that $A(r)$ is finite for all choices of $\alpha\geq
0,\beta>0,$ whereas $B(r)$ is finite if and only if $\alpha<p-1.$
Both functions $A,B$ are locally bounded and continuous next to
$0$ and $\infty$. Moreover, for $r$ close to $0$, we have
\[
A(r)^{\frac{1}{p}}(B(r))^{1-\frac{1}{p}}\asymp \left(
r^{-1+\frac{\alpha}{p}+ \frac{1}{p}}+C\right)\cdot \left(
r^{-\frac{\alpha}{p} +1-\frac{1}{p}}\right) < Const.
\]
Therefore (\ref{mucexp}) holds true if and only if
$\ds\limsup_{r\to\infty} A(r)B(r)^{p-1}<\infty.$

By a direct application of the  de l'Hospital rule we see that for
$a\in\mathbb{R}, b>0$
\[
\lim_{r\to\infty}\frac{\int_r^\infty x^a
e^{-{x^b}}dx}{r^{a+1-b}e^{-r^b}}=\lim_{r\to\infty}
\frac{r^ae^{-r^b}}{br^ae^{-r^b}-(a+1-b)r^{a-b}e^{-{r^b}}}=\frac{1}{b},
\]
and for $a<1, C>0$
\[
\lim_{r\to\infty}\frac{\int_0^r \frac{e^{Cx^b}}{x^a}
dx}{e^{Cr^b}r^{-(a+b-1)}}=\lim_{r\to\infty}\frac{e^{Cr^b}}{r^a}\left[
\frac{bCe^{Cr^b}}{r^a}-\frac{(a+b-1)e^{Cr^b}}{r^{a+b}}\right]^{-1}=\frac{1}{bC}.
\]
Therefore, for $r$ large, we have
\begin{eqnarray*}
A(r)\asymp r^{(\beta-1)(p-1)+\alpha}e^{-r^\beta}&\mbox{ and }&
B(r)\asymp
r^{-\frac{\alpha}{p-1}-(\beta-1)}e^{\frac{r^\beta}{p-1}},
\end{eqnarray*}
and so $A(r)B(r)^{p-1}\asymp Const$ for large values of $r.$
Therefore (\ref{mucexp}) is fulfilled whenever $0\leq \alpha<p-1,
\beta>0. $

We end up with the following theorem.

\begin{theo}\label{Mukho1}
Let $p>1$ and let $\mu(dx)= x^{\alpha}e^{-x^{\beta}}dx$, where
$\alpha\neq p-1,\,\beta=0,$ or  $0\leq \alpha<p-1, \beta>0.$
 Suppose that $M$ is a an $N-$function
satisfying condition {\bf (M)}, and  $Q$ is another $N-$function,
such that
\[
\frac{M(t)}{t^2}rs\le M(t) +c r^p +Q(s),\ \hbox{\rm for every}\
t,r,s>0.
\]
Then inequalities (\ref{iii}) and (\ref{iiin}) hold for any $u\in
C_0^\infty(\rp),$ with constants $K_1,K_2,K$ independent of $u.$
\end{theo}
\end{ex}

\begin{rem}\rm
As to  the validity of Orlicz-space counterparts of (\ref{inek1}),
 which would then yield {\bf (H)}, we refer to the
papers of Bloom-Kerman \cite{bloker1,bloker}, Lai \cite{lai4},
Heinig-Maligranda \cite{hm}, Bloom-Kerman \cite{bloker},
Heinig-Lai \cite{lai6}, their references and also to the authors'
paper \cite{akikppstudia}, Section 3.3, where another type of
sufficient conditions for {\bf (H)} to hold on $\r_+$ is given.
\end{rem}

\subsubsection{Multidimensional Hardy inequalities}

{\bf Inequalities on bounded domains}

\smallskip

\noindent The multidimensional Hardy inequalities of the form
\begin{equation}\label{harta}
\int_{\Omega} \left(\frac{1}{\delta (x)} |u(x)|
\right)^q\delta(x)^a dx\le C \int_{\Omega} |\nabla
u(x)|^q\delta(x)^adx,\ \ u\in C_0^1(\Omega),
\end{equation}
where $\Omega\subseteq\rn$ is a bounded domain with sufficiently
regular boundary,  $a<q-1$, $1<q<\infty$, $\delta(x)={\rm
dist}(x,\partial\Omega)$, were first obtained by Necas \cite{ne}
(for bounded domains with Lipschitz boundary) and extended further
by Kufner (\cite{kufner}, Theorem 8.4) and Wannebo \cite{wan} to
H\"older domains.

 As a direct consequence
of of Theorem \ref{theoA}, we obtain Gagliargo-Nirenberg-type
inequalities within $L^p$-spaces with the distance from the
boundary, which can be stated as follows.

\begin{theo}\label{distan} Suppose that $\Omega$ is a bounded
Lipschitz domain and $q>1,$ $a>q-1$ . Then we have
\begin{description}
\item[i)]
if $p\geq 2,$ $r>1,$ $\frac{2}{p}=\frac{1}{q}+\frac{1}{r}$ then
for every $u\in C_0^\infty(\Omega)$ one has:
\begin{equation}\label{zabka}
\left( \int_{\Omega} |\nabla
u(x)|^p\delta(x)^adx\right)^{\frac{2}{p}}\le c  \left(
\int_{\Omega} |\nabla^{(2)}
u(x)|^q\delta(x)^adx\right)^{\frac{1}{q}}\cdot \left(
\int_{\Omega} | u(x)|^r\delta(x)^adx\right)^{\frac{1}{r}},
\end{equation}
with a constant $c>0$ independent of $u.$
\item[ii)] if $M$ and $Q$ are
    $N-$functions such that
    \[\frac{M(u)}{u^2}\,vw\leq M(u)+cv^q +Q(w),\]
    and $M$ satisfies condition {\bf (M)}, then
for every $u\in C_0^\infty(\Omega)$ one has:
    \begin{equation}\label{orli}
\left( \int_{\Omega} M(|\nabla u(x)|)\delta(x)^adx\right)\le c
\left( \int_{\Omega} |\nabla^{(2)} u(x)|^q\delta(x)^a dx +
\int_{\Omega} Q(| u(x)|) \delta(x)^adx\right),
\end{equation}
and \begin{equation}\label{orli1} \|\nabla
u\|^2_{L^M(\Omega,\mu)}\leq
c\|\nabla^{(2)}u\|_{L^q(\Omega,\mu)}\|u\|_{L^Q(\Omega,\mu)},
\end{equation}
with constants independent on $u$.
\end{description}
\end{theo}

{\bf Proof.} {\bf i)} Obviously, $M(\lambda)=\lambda^p$ satisfies
condition {\bf (M)} and $P(\lambda)=\lambda^q,
Q(\lambda)=\lambda^r$ satisfy {\bf (Y)} due to Remark
\ref{konsul}. Moreover, the measure
 $\mu (dx) ={\rm exp}(-\varphi (x))dx$ where
 $\varphi (x)=-a \ln \delta (x)$ satisfies {\bf ($\mu$)}.
 It is not hard to verify that $|\nabla \varphi (x)|\asymp \frac{1}{\delta (x)}$
 (note that $|\nabla\delta|\asymp Const$). This together with (\ref{harta})
  implies {\bf (H)}. Now it suffices to apply Theorem \ref{theoA}.
  Part {\bf ii)} is proven similarly.
  \hfill$\Box$

\begin{rem}\rm
Note that $\delta (x)^a$ is an $A_q$-weight when $-1<a<q-1$ (see
e.g. \cite{tor}  for definition of $A_p$ weights introduced by
Muckenhoupt). Gagliardo-Nirenberg inequalities with $A_p$-weights
within
    homogeneous spaces were earlier obtained in \cite{skch94}, \cite{akikppstud} by  different methods.
    Those results also covered the case $1<p<2.$
\end{rem}

\begin{rem} \rm Counterparts of inequality  (\ref{harta}) in Orlicz norms
were obtained by Cianchi in  \cite{cia1}.
\end{rem}

\bigskip
 \noindent {\bf Inequalities on
$\rn$}

\noindent Hardy inequality on $\rn$ with power weights (see e.g.
\cite{ckn}, \cite{ma}, page 70 in \cite{kmp}, and their
references)
\[
\| |x|^{\alpha -1}|u|\|_{L^q}\le C\| |x|^\alpha |\nabla
u|\|_{L^q},
\]
where $u\in C_0^\infty (\rn)$, $\frac{1}{q}+\frac{\alpha
-1}{n}>0$, $q>1$, give rise to the Gagliardo-Nirenberg
inequalities on $\rn$  with power weights $|x|^\alpha$ and
$N-$functions $M,P=P(\lambda )=\lambda^q,Q,$  satisfying {\bf
(Y)}. The result, obtained directly from Theorem \ref{theoA} reads
as follows.

\begin{theo}\label{distances} Suppose that $\frac{1}{q}+\frac{\alpha
-1}{n}>0$, $q>1$. Then we have
\begin{description}
\item[i)]
if $p\geq 2,$ $r>1,$ $\frac{2}{p}=\frac{1}{q}+\frac{1}{r}$ then
for every $u\in C_0^\infty(\Omega)$ one has:
\begin{equation}\label{zabka1}
\left( \int_{\Omega} |\nabla
u(x)|^p|x|^{\alpha}dx\right)^{\frac{2}{p}}\le c  \left(
\int_{\Omega} |\nabla^{(2)} u(x)|^q|x|^{\alpha}
dx\right)^{\frac{1}{q}}\cdot \left( \int_{\Omega} |
u(x)|^r|x|^{\alpha}dx\right)^{\frac{1}{r}},
\end{equation}
with a constant $c>0$ independent of $u.$
\item[ii)] if $M$ and $Q$ are
    $N-$functions such that
    \[\frac{M(u)}{u^2}\,vw\leq M(u)+cv^q +Q(w),\]
    and $M$ satisfies condition {\bf (M)}, then
for every $u\in C_0^\infty(\Omega)$ one has:
    \begin{equation}\label{orlii}
\left( \int_{\Omega} M(|\nabla u(x)|)|x|^{\alpha}dx\right)\le c
\left( \int_{\Omega} |\nabla^{(2)} u(x)|^q|x|^{\alpha} dx +
\int_{\Omega} Q(| u(x)|) |x|^{\alpha}dx\right),
\end{equation}
and \begin{equation}\label{orlii1} \|\nabla
u\|^2_{L^M(\Omega,\mu)}\leq
c\|\nabla^{(2)}u\|_{L^q(\Omega,\mu)}\|u\|_{L^Q(\Omega,\mu)},
\end{equation}
with constants independent on $u$.
\end{description}
\end{theo}

\subsection{The scope of Theorem \ref{theoAB}}

We are now to discuss the validity of Theorem \ref{theoAB}. To
abbreviate the discussion we only focus on the condition {\bf
(H1)} in its assumptions. Contrarily to the Hardy inequality {\bf
(H)} the condition  {\bf (H1)} barely seems to be investigated.

\subsubsection{Inequalities on bounded domains}\label{io}

 {\bf (A) Result by Oinarov.} Condition {\bf(H1)} and its special
variant (\ref{assumB1}) has drawn less attention than the
`classical' Hardy inequality {\bf (H)}. In his paper  \cite{oi},
Oinarov considered inequalities
\begin{equation}\label{oina}
\left( \int_a^b |\omega  u|^q dr\right)^{\frac{1}{q}}\le C\left(
\int_a^b |\rho u'|^p dr + \int_a^b |vu|^p dr
\right)^{\frac{1}{p}},
\end{equation}
for general weights $\omega,v, \rho$ and derived necessary and
sufficient conditions needed for (\ref{oina}) to hold for all
$u\in C_0^\infty(a,b).$  Our condition (\ref{assumB1}) with
$\mu(dx)={\rm e}^{-\vp(x)}dx$ is exactly (\ref{oina}), under the
choice of $\omega(r)= \varphi'(r){\rm e}^{-\frac{\varphi(r)}{p}},$
$v(r)=\rho(r)={\rm e}^{-\frac{\varphi(r)}{p}}, p=q.$

\smallskip

\noindent {\bf (B) Result by Cianchi.} Orlicz-norms counterpart of
{\bf (H1)}:
\begin{equation}\label{chichi}
\| \frac{u}{d^{1+\alpha}}\|_{L^B(G)}\le C\left( \|
\frac{u}{d^{\alpha}}\|_{L^A(G)}+\| \frac{\nabla
u}{d^{\alpha}}\|_{L^A(G)} \right)
\end{equation}
 has been established by Cianchi in \cite{cia1}.
Here $G\subset \rn$ is a sufficiently regular domain, $A$ and $B$
are $N-$functions related by a certain domination condition (in
particular it it possible to have $A=B$),
$d(x)=\mbox{dist}\,(x,\partial G)$, and the measure considered is
the Lebesgue measure. In our work, we
 need modular versions of (\ref{chichi}); in general they do not
come as its direct consequence. Note that in the case of
homogeneous $N$-functions (\ref{chichi}) is an extension of
(\ref{harta}).

\smallskip

\noindent {\bf (C) Authors' approach. Inequalities on intervals.}
In the forthcoming paper \cite{akikpp09}, the authors  work out
condition concerning $M,\varphi$ which ensure the validity of
(\ref{assumB1}) on intervals $(a,b)\subset \r,$ not excluding the
cases $a=-\infty$ or $b=\infty.$  See also the paper \cite{KOKPP}
devoted solely to the Hardy and Gagliardo-Nirenberg inequalities
for Gaussian measure on $\rn.$

\noindent Email addresses: {\tt kalamajs@mimuw.edu.pl}\\
\phantom{ } \hskip 2.75cm  {\tt kpp@mimuw.edu.pl}

\end{document}